\documentclass[12pt,leqno]{article}
\usepackage{amsmath,amssymb,amsthm,amscd}
\usepackage{dsfont}
\numberwithin{equation}{section} \allowdisplaybreaks
\begin{document}
\newtheorem{theorem}{Theorem}[section]
\newtheorem{defin}{Definition}[section]
\newtheorem{prop}{Proposition}[section]
\newtheorem{corol}{Corollary}[section]
\newtheorem{lemma}{Lemma}[section]
\newtheorem{rem}{Remark}[section]
\newtheorem{example}{Example}[section]
%\label{} %\ref{}
\title{Dirac Structures and Generalized Complex Structures on
$TM\times\mathds{R}^h$}
\author{{\small by}\vspace{2mm}\\Izu Vaisman}
\date{}
\maketitle
{\def\thefootnote{*}\footnotetext[1]%
{{\it 2000 Mathematics Subject Classification: 53C15, 53D17}. 
}}
\begin{center} \begin{minipage}{12cm}
A{\footnotesize BSTRACT. We consider Courant and Courant-Jacobi
brackets on the stable tangent bundle $TM\times\mathds{R}^h$ of a
differentiable manifold and corresponding Dirac, Dirac-Jacobi and
generalized complex structures. We prove that Dirac and
Dirac-Jacobi structures on $TM\times\mathds{R}^h$ can be prolonged
to $TM\times\mathds{R}^k$, $k>h$, by means of commuting
infinitesimal automorphisms. Some of the stable, generalized,
complex structures are a natural generalization of the normal,
almost contact structures; they are expressible by a system of
tensors $(P,\theta,F,Z_a,\xi^a)$ $(a=1,...,h)$, where $P$ is a
bivector field, $\theta$ is a $2$-form, $F$ is a $(1,1)$-tensor field,
$Z_a$ are vector fields and $\xi^a$ are $1$-forms, which satisfy
conditions that generalize the conditions satisfied by a normal,
almost contact structure $(F,Z,\xi)$. We prove that such a generalized
structure projects to a generalized, complex structure of a space
of leaves and characterize the structure by means of the projected
structure and of a normal bundle of the foliation. Like in the
Boothby-Wang theorem about contact manifolds, principal torus
bundles with a connection over a generalized, complex manifold
provide examples of this kind of generalized, normal, almost
contact structures.}
\end{minipage}
\end{center} \vspace{5mm}
\section{Brackets on stable tangent bundles}
All the manifolds and mappings of the present note are assumed of
the $C^\infty$ class. In differential topology, the vector bundle
$TM\times\mathds{R}^h$, where $M$ is a differentiable manifold, is
called the {\it stable tangent bundle}. The name comes from the
fact that $h$ may be arbitrary and the main interest is for
objects defined up to an equivalence that involves a change of $h$
(e.g.,
\cite{GGK}, Appendix D). On the other hand, interesting
differential-geometric structures of $M$ may be defined via a
stable tangent bundle; for instance, an almost contact structure
is equivalent with an almost complex structure on
$TM\times\mathds{R}$ and the integrability of the latter
characterizes the normality of the former
\cite{Bl1}. The aim of the present paper is to extend the notions
of a Dirac structure
\cite{{CW},{C}}, a Dirac-Jacobi structure
\cite{{Wd},{GM}} and a generalized complex structure \cite{{H},{G}}
to $TM\times\mathds{R}^h$. We will give some results concerning
the prolongation of a structure from $TM\times\mathds{R}^h$ to
$TM\times\mathds{R}^k$, $k>h$ and discuss structures defined on $M$  
by generalized complex structures in $TM\times\mathds{R}^h$.

Our general notation is $\chi^k(M)$ for the space of $k$-vector
fields, $\Omega^k(M)$ for the space of differential $k$-forms,
$\Gamma$ for the space of global cross sections of a vector
bundle, $X,Y,..$ for contravariant vectors (or vector fields),
$\alpha,\beta,...$ for covariant vectors (or $1$-forms), $u,v,...$
for vectors of $\mathds{R}^h$ and a dot for the natural scalar
product of $\mathds{R}^h$. Furthermore, we define the {\it big
tangent bundle} $T^{big}M=TM\oplus T^*M$ and the {\it stable big
tangent bundle of index $h$}
\begin{equation}\label{bigstable} \mathbf{T}_{h}^{big}M=
(TM\times\mathds{R}^h)\oplus(T^*M\times\mathds{R}^h),
\end{equation} where $h$ is any non negative integer,
which, usually, we shall not mention (on the other hand, notice
the boldface character, except for $\mathbf{T}_{0}^{big}M=T^{big}M$).

The bundle $\mathbf{T}^{big}M$ has a natural, neutral metric
(i.e., non degenerate and of signature zero)
\begin{equation}\label{gbig}
G((X_1,u_1)\oplus(\alpha_1,v_1),(X_2,u_2)\oplus(\alpha_2,v_2))\end{equation}
$$=\frac{1}{2}(\alpha_2(X_1)+\alpha_1(X_2)+u_1\cdot v_2+u_2\cdot
v_1),$$ where $\oplus$ denotes the addition of elements in
different terms of a direct sum of vector spaces. Similarly, and with the generic
notation
\begin{equation}\label{Xrond} \mathcal{X}=
(X,u)\oplus(\alpha,v),
\end{equation} one has the non degenerate $2$-form
\begin{equation}\label{omegabig}
\Omega( \mathcal{X}_1,\mathcal{X}_2)
=\frac{1}{2}(\alpha_2(X_1)-\alpha_1(X_2)+u_1\cdot v_2-u_2\cdot
v_1).\end{equation}

The usual, twisted, Courant bracket is an operation on
$\Gamma(T^{big}M)$ (i.e., $h=0$) defined by \cite{{CW},{C},{SW}}
\begin{equation}\label{crosetCSW}
[X\oplus\alpha,Y\oplus\beta]_\Phi=[X,Y]\oplus[L_X\beta-L_Y\alpha
\end{equation}
$$-d(\Omega(X\oplus\alpha,Y\oplus\beta))+i(X\wedge Y)\Phi],$$
where $L$ denotes the Lie derivative and $\Phi$ is a closed $3$-form; if $\Phi=0$ the bracket is untwisted.

We shall also need the notion of a {\it conformal change} of the
Courant bracket, which is a particular case of the deformation by a
$1$-cocycle defined by Grabowski and Marmo \cite{GM} and a
generalization of a change used by Wade \cite{Wd2} and by Petalidou and Nunes da Costa
\cite{PC}. An automorphism $ \mathcal{C}_{\tau}: T^{big}M
\rightarrow T^{big}M$ defined by
\begin{equation}\label{confchange} \mathcal{C}_{\tau}
(X\oplus\alpha)= X\oplus e^\tau\alpha,\hspace{3mm}
\tau\in C^\infty(M),\end{equation}
will be called a conformal change of $T^{big}M$, because it
produces a conformal change of the metric $G$ in the sense that
\begin{equation}\label{confpeg} G(\mathcal{C}_{\tau}
(X\oplus\alpha), \mathcal{C}_{\tau} (Y\oplus\beta)) = e^{\tau}
G(X\oplus\alpha,Y\oplus\beta). \end{equation} Correspondingly, the
bracket
\begin{equation}\label{confgen}
[X\oplus\alpha,Y\oplus\beta]_{\tau,\Phi} =
\mathcal{C}_{-\tau}[\mathcal{C}_{\tau}(X\oplus\alpha),
\mathcal{C}_{\tau}(Y\oplus\beta)]_\Phi,
\end{equation}
where the right hand side is  defined by means of
(\ref{crosetCSW}), will be called a {\it conformal-Courant
bracket}. Using (\ref{crosetCSW}) we get
\begin{equation}\label{conf2}
[X\oplus\alpha,Y\oplus\beta]_{\tau,\Phi} =
[X\oplus\alpha,Y\oplus\beta]_{e^{-\tau}\Phi}
\end{equation} $$+(0\oplus[(X\tau)\beta-(Y\tau)\alpha
-\Omega(X\oplus\alpha,Y\oplus\beta)d\tau]).$$
\begin{rem}\label{nonC} {\rm From
(\ref{confgen}), it follows that the conformal Courant bracket
satisfies the following among the axioms of a Courant algebroid of
anchor $\rho=pr_{TM}$ \cite{LWX}:\\

i) $\rho[X\oplus\alpha,Y\oplus\beta]_{\tau,\Phi}=
[\rho(X\oplus\alpha),\rho(Y\oplus\beta)]$,\\

ii) $\rho(\partial f)=0\hspace{3mm}(f\in C^\infty(M), \partial
f=(0,df))$,\\

iii) $[X\oplus\alpha,f(Y\oplus\beta)]_{\tau,\Phi} =
f[X\oplus\alpha,Y\oplus\beta]_{\tau,\Phi}+(Xf)(Y\oplus\beta)$\vspace{2mm}\\ 
\hspace*{5cm}$-G(X\oplus\alpha,Y\oplus\beta)\partial f.$\vspace{2mm}\\ But, the
other axioms of Courant algebroids are not satisfied. In
particular, formula (\ref{confgen}) yields
\begin{equation}\label{JnonC} \sum_{Cycl(1,2,3)}
[[X_1\oplus\alpha_1,X_2\oplus\alpha_2]_{\tau,\Phi},
X_3\oplus\alpha_3]_{\tau,\Phi} =\frac{1}{3}\partial\sum_{Cycl(1,2,3)}
G([X_1\oplus\alpha_1,\end{equation} $$X_2\oplus\alpha_2]_{\tau,\Phi},
X_3\oplus\alpha_3) +\frac{1}{3}\sum_{Cycl(1,2,3)}
G([X_1\oplus\alpha_1,X_2\oplus\alpha_2]_{\tau,\Phi}, X_3\oplus\alpha_3)
\partial\tau.$$}\end{rem}

We define bracket operations on the stable big tangent bundle from
the brackets (\ref{crosetCSW}) and (\ref{confgen}) on the manifold
$M\times\mathds{R}^h$. First, a vector field $\tilde
X\in\chi^1(M\times\mathds{R}^h)$, a $1$-form
$\tilde\alpha\in\Omega^1(M\times\mathds{R}^h)$, etc. will be
called {\it translation invariant} if they are preserved by the
natural action of the group of translations of $
\mathds{R}^h$ on $M\times\mathds{R}^h$, equivalently, if
\begin{equation}\label{invarsect} L_{\frac{\partial}{\partial
t^a}}\tilde X=0, L_{\frac{\partial}{\partial t^a}}\tilde\alpha=0,
\hspace{3mm}{\rm etc.},
\end{equation} where $t^a$ $(a=1,...,h)$ are
the natural coordinates of $ \mathds{R}^h$. Obviously, the space
$\Gamma\mathbf{T}_{h}^{big}M$ is naturally isomorphic with the
space of translation invariant cross sections of
$T^{big}(M\times\mathds{R}^h)$ by the identification
\begin{equation}\label{isomtau2} (X,u)\oplus(\alpha,v)\leftrightarrow
(X+\sum_{a=1}^hu^a\frac{\partial}{\partial t^a})\oplus
(\alpha+\sum_{a=1}^hv_adt^a).\end{equation}

The restriction of the Courant bracket (\ref{crosetCSW}) of
$M\times\mathds{R}^h$ with a twist form
\begin{equation}\label{twistform}
\Phi+\sum_{a=1}^hdt^a\wedge\Psi^a,\hspace{5mm}
\Phi\in\Omega^3(M),\Psi^a\in\Omega^2(M),\,
d\Phi=0,d\Psi^a=0\end{equation} to translation invariant cross
sections, composed by (\ref{isomtau2}), will be called the {\it
stable Courant bracket of index $h$} on $M$.

It follows that the stable Courant bracket is given by
\begin{equation}\label{crosettwbig}
[(X_1,u_1)\oplus(\alpha_1,v_1),(X_2,u_2)\oplus(\alpha_2,v_2)]_C
=([X_1,X_2],X_1u_2-X_2u_1) \end{equation}
$$\oplus(L_{X_1}\alpha_2-L_{X_2}\alpha_1+
\frac{1}{2}d(\alpha_1(X_2)-\alpha_2(X_1))+i(X_1\wedge
X_2)\Phi+u_1\cdot (i(X_2)\Psi)$$
$$-u_2\cdot (i(X_1)\Psi)
+\frac{1}{2}(u_2\cdot dv_1+v_2\cdot du_1-u_1\cdot dv_2-v_1\cdot
du_2),X_1v_2-X_2v_1+\Psi(X_1,X_2)),$$ where $\Psi$ is the $\mathds{R}^h$-valued form of components $\Psi^a$. The index $C$ will be
omitted if there is no danger of confusion.
\begin{rem}\label{obsalgtransitiv} {\rm The bundle
$\mathbf{T}_{h}^{big}M$ with the metric $G$, the bracket
(\ref{crosettwbig}) and the anchor $\rho=pr_{TM}$ is a transitive
Courant algebroid. Thus, the bracket (\ref{crosettwbig}) could
have been derived from the general formulas of \cite{VC8}.}
\end{rem}

Furthermore, if we replace the Courant bracket by the conformally
changed bracket (\ref{confgen}) of $M\times\mathds{R}^h$ with the
change function $\tau=\tau(t^1,...,t^h)$, restricted to
translation invariant sections and with the result taken at
$t^a=0$, we get a new bracket, which has the following expression
\begin{equation}\label{crosetwtbig}
[(X_1,u_1)\oplus(\alpha_1,v_1),(X_2,u_2)\oplus(\alpha_2,v_2)]_W
=([X_1,X_2],X_1u_2-X_2u_1) \end{equation}
$$\oplus(L_{X_1}\alpha_2-L_{X_2}\alpha_1+
\frac{1}{2}d(\alpha_1(X_2)-\alpha_2(X_1))+d_0\tau(u_1)\alpha_2
-d_0\tau(u_2)\alpha_1$$
$$+e^{-\tau(0)}i(X_1\wedge X_2)\Phi+e^{-\tau(0)}u_1\cdot (i(X_2)\Psi)
-e^{-\tau(0)}u_2\cdot (i(X_1)\Psi)$$ $$ +\frac{1}{2}(u_2\cdot
dv_1+v_2\cdot du_1-u_1\cdot dv_2-v_1\cdot du_2),
X_1v_2-X_2v_1+e^{-\tau(0)}\Psi(X_1,X_2)$$ $$
+\frac{1}{2}(\alpha_1(X_2)-\alpha_2(X_1)+u_2\cdot v_1-u_1\cdot
v_2)d_0\tau+d_0\tau(u_1)v_2-d_0\tau(u_2)v_1).$$

The bracket (\ref{crosetwtbig}) will be called the {\it Wade (Courant-Jacobi)
bracket of index $h$} because it was first discovered by Wade
\cite{Wd} in the case $h=1,\tau=t$. The index $W$ will be omitted if
there is no danger of confusion.
\begin{rem}\label{obsuzero} {\rm In the
untwisted case, and if $\partial\tau/\partial t^a=const.$, the
bracket (\ref{crosetwtbig}) with an arbitrary $t$ is translation
invariant. Hence, in this case the translation invariant section
defined by
$[(X_1,u_1)\oplus(\alpha_1,v_1),(X_2,u_2)\oplus(\alpha_2,v_2)]_W$
is equal to the conformal-Courant bracket
$$
[(X_1+\sum_{a=1}^hu^a_1\frac{\partial}{\partial t^a})
\oplus(\alpha_1+\sum_{a=1}^hv_{1,a}dt^a),
(X_2+\sum_{a=1}^hu^a_2\frac{\partial}{\partial t^a})
\oplus(\alpha_2+\sum_{a=1}^hv_{2,a}dt^a]_\tau. $$}
\end{rem}
\section{Stable Dirac and Dirac-Jacobi Structures}
A twisted, respectively untwisted, {\it stable Dirac structure of
index $h$} is a maximal, $G$-isotropic subbundle $
L\subseteq\mathbf{T}_h^{big}M$, which satisfies the integrability
condition of being closed by the twisted, respectively untwisted,
bracket (\ref{crosettwbig}). If the integrability condition is not
satisfied the structure is {\it stable almost Dirac}. In the case
$h=0$ we have the usual Dirac structures of \cite{{CW},{C}}.

We can define an equivalence relation that justifies the ``stable"
terminology. A stable almost Dirac structure of index $h$,
$L\subseteq\mathbf{T}_h^{big}M$, extends to the following
structures of index $h+k$:
\begin{equation}\label{trivialext2}\begin{array}{l}\hat L_1=
\{(X,u,w)\oplus(\alpha,v,0)\,/\,
(X,u)\oplus(\alpha,v)\in L,\,w\in\mathds{R}^k\},\vspace{2mm}\\
\hat L_2=\{(X,u,0)\oplus(\alpha,v,w)\,/\, (X,u)\oplus(\alpha,v)\in
L,\,w\in\mathds{R}^k\}.\end{array}\end{equation} If the $\mathds{R}^h$-valued $2$-form $\Psi$ of (\ref{twistform}) is extended by $k$ zero components to a $\mathds{R}^{h+k}$-valued form, the structures $L,\hat L_1,\hat L_2$ simultaneously are integrable or not. Notice also the existence of the metric preserving automorphism $F: 
\mathbf{T}_{h+k}^{big}M=\mathbf{T}_h^{big}M\oplus\mathds{R}^{2k}
\rightarrow \mathbf{T}_{h+k}^{big}M$ defined by $F|_{\mathbf{T}_h^{big}M}
=Id,F|_{\mathds{R}^{2k}}(w,0)=(0,w)$, which sends $\hat L_1$ onto $\hat L_2$ and conversely.
\begin{defin}\label{echivalenta} {\rm The
(almost) Dirac structures $L\subseteq\mathbf{T}_h^{big}M$,
$L'\subseteq\mathbf{T}_{h'}^{big}M$ are called {\it stably equivalent}
if there are non negative integers $k,k'$ such that $h+k=h'+k'$, and
there exists a metric preserving, bundle automorphism
$\varphi:(\mathbf{T}_{h+k}^{big}M,G) \rightarrow
(\mathbf{T}_{h'+k'}^{big}M,G)$ that sends the prolongation $\hat
L_1$ defined by (\ref{trivialext2}) to the similar prolongation
$\hat L'_1$.}\end{defin}

It is trivial to see that stable equivalence is an equivalence
relation. Moreover, if $\varphi$ exists then $F'\circ\varphi\circ
F$, where $F$ was defined above and $F'$ is defined similarly, sends $\hat L_2$ onto $\hat L'_2$ and conversely.
\begin{example}\label{exhsubvar} {\rm 
Let $M^m$ be a differentiable manifold
and $N^n$ be a submanifold with a trivial normal bundle $T_NM/TN$.
Then, there exists a non canonical isomorphism
\begin{equation}\label{isomptred} T_NM\approx TN\oplus
\mathds{R}^{m-n}.\end{equation} The isomorphism (\ref{isomptred})
yields a $G$-preserving isomorphism
\begin{equation}\label{isoinexemplu} I:
T^{big}_NM\stackrel{\approx}{\rightarrow}\mathbf{T}^{big}_{m-n}N,\end{equation} and if $D$ is an almost Dirac structure of
$M$ then $I(D|_N)$ is a stable almost Dirac structure of $N$. Since the
isomorphism (\ref{isoinexemplu}) is not unique, it is rather the corresponding equivalence class of stable almost Dirac structures that is well defined. Generally, the
integrability of $D$ does not imply the integrability of $I(D|_N)$.}\end{example}

The meaning of the notion of a stable Dirac structure of index $h$
is given by the following simple proposition.
\begin{prop}\label{prophogen}
A stable, almost Dirac structures $L$ of $M$ may be identified
with a translation invariant, almost Dirac structure
$\tilde{L}$of the manifold $M\times\mathds{R}^h$. The structures
$L$ and $\tilde{L}$ are simultaneously integrable or
not.\end{prop}
\begin{proof} The invariance of $\tilde{L}$
by translations means that
$\forall s\in\mathds{R}^h,\forall(x,t)\in M\times\mathds{R}^h$ the translation $\tau_s(x,t)=(x,t+s)$ satisfies the condition
$$\tilde{L}_{(x,t+s)}=(\tau_s)_*(\tilde L_{(x,t)})=
\{(\tau_s)_{*(x,t)}\tilde X\oplus(\tau_{-s})^*_{(x,t)}\tilde\alpha\,/\,
\tilde X\oplus\tilde\alpha\in\tilde L_{(x,t)}\}.$$
This condition is equivalent with
$$ L_{\frac{\partial}{\partial
t^a}}\tilde{\mathcal{X}}\in\Gamma\tilde{L},\hspace{5mm}\forall \tilde{\mathcal{X}}\in\Gamma\tilde{L},
$$ where the Lie derivative is applied to each
component of $\tilde{\mathcal{X}}$. Furthermore, the translation invariance condition is also equivalent with the fact that $\Gamma\tilde L$ has local bases that consist of translation invariant cross sections $Z_i\oplus\theta_i$ of the form (\ref{isomtau2}); these bases are obtained by translating local bases of cross sections of $\tilde L|_{t=0}$.

Now, the stable almost Dirac structure
$L\subseteq\mathbf{T}^{big}M$ produces the translation invariant almost Dirac structure
\begin{equation}\label{deftildeL}
\tilde L_{(x,t)}= \{(\tau_t)_{*(x,0)}X\oplus(\tau_{-t})^*_{(x,0)}
\alpha\,/\,X\oplus\alpha\in L_x\}\end{equation} of 
$M\times \mathds{R}^h$.
Moreover, all the translation invariant almost Dirac structures of $M\times
\mathds{R}^h$ are produced in this way by $L=\tilde L|_{t=0}$; this may be seen by using
the local invariant bases of a translation invariant structure. Furthermore, $\Gamma L$ is identifiable with the space of
the translation invariant cross sections of $\tilde L$ of (\ref{deftildeL}) via
(\ref{isomtau2}). Hence, the integrability of $\tilde{L}$ implies the
integrability of $L$. The converse follows by expressing the cross
sections of $\tilde{L}$ by means of local, translation invariant
bases and by using the properties of the Courant bracket and the
isotropy of $\tilde{L}$.
\end{proof}
\begin{rem}\label{cazneomogen} {\rm If $\tilde{L}$ is an arbitrary
almost Dirac structure of $M\times\mathds{R}^h$,
$L=\tilde{L}|_{t=0}$ is a stable almost Dirac structure of $M$,
which may not be integrable even if $\tilde{L}$ is integrable
because the restriction to $t=0$ of the Courant bracket on
$M\times\mathds{R}^h$ is not the stable Courant bracket on $M$,
generally.}\end{rem}
\begin{example}\label{exemlulPi} {\rm Let
\begin{equation}\label{Pi} \Pi=
W+\sum_{a=1}^hV_a\wedge\frac{\partial}{\partial t^a},
\hspace{2mm}W\in\chi^2(M),V_a\in\chi^1(M)
\end{equation} be a Poisson structure on $M\times
\mathds{R}^h$ with the twist form (\ref{twistform}).
The Poisson condition $[\Pi,\Pi]=0$
(Schouten-Nijenhuis bracket) is equivalent with the conditions
\begin{equation}\label{condptPi} [W,W]=0,\hspace{2mm}L_{V_a}W=0,
\hspace{2mm} [V_a,V_b]=0\hspace{5mm}(a,b=1,...,h), \end{equation}
i.e., $W$ is a Poisson structure on $M$ and $V_a$ are $h$
commuting, infinitesimal automorphisms of $W$. The graph of $\sharp_\Pi$
($\sharp_P$ is defined by $\beta(\sharp_\Pi\alpha)=\Pi(\alpha,\beta)$) is a translation
invariant Dirac structure $\tilde L$ on $M\times \mathds{R}^h$ obtained by the translation of the stable Dirac structure
\begin{equation}\label{Linex1} L=\tilde{L}|_{t=0}=\{(\sharp_W\alpha-u\cdot
V,\alpha(V))\oplus(\alpha,u)\},
\end{equation} where $\alpha\in T^*M,u\in \mathds{R}^h$,
$V$ is the $ \mathds{R}^h$-valued vector field on $M$ with the
components $V_a$ (i.e., $V_xf=(V_{a,x}f)\in\mathds{R}^h,\,\forall f\in C^\infty(M),\,\forall x\in M$) and $\alpha(V)\in\mathds{R}^h$ is calculated on
the components $V_a$ of $V$.}\end{example}
\begin{example}\label{exemplulTheta} {\rm
The graph of $\flat_\Theta$ $(\flat_\Theta(X)=i(X)\Theta)$ where
\begin{equation}\label{Theta}
\Theta=\sigma+\sum_{a=1}^h\theta_a\wedge dt^a\;\;(\sigma\in\Omega^2(M),
\theta_a\in\Omega^1(M),d\sigma=0,
d\theta_a=0)\end{equation} is a translation invariant Dirac
structure $\tilde L$ with a corresponding stable, Dirac structure
\begin{equation}\label{LptTheta} L=\{(X,u)\oplus(\flat_\sigma
X-u\cdot\theta,\theta(X))\} \end{equation} on $M$.} \end{example}

The following prolongation theorem extends Example
\ref{exemlulPi}.
\begin{theorem}\label{propex}
Let $L\subseteq \mathbf{T}_h^{big}M$ be an untwisted, stable Dirac
structure of index $h$ on $M$ and let $V_p$ $(p=1,...,k)$ be
commuting infinitesimal automorphisms of $L$. Then, formula
\begin{equation}\label{Linex2} \hat L_x=\{(X-w\cdot V(x),u,\alpha(V(x)))
\oplus(\alpha,v,w)\end{equation}
$$\,/\,(X,u)\oplus(\alpha,v)\in L,\,w\in\mathds{R}^k\} \hspace{2mm}(x\in M)$$ defines an untwisted,
stable Dirac structure $\hat L$ of index $h+k$.\end{theorem}
\begin{proof} Notice that the elements and cross sections of $\hat L$ may be written as
\begin{equation}\label{descinex2}
[(X,u,\alpha(V))\oplus(\alpha,v,0)] + [(-w\cdot
V,0,0)\oplus(0,0,w)],\end{equation} where $(X,u)\oplus(\alpha,v)$ are either elements or cross sections of $L$ and $w\in\mathds{R}^k$, $w\in C^\infty(M,\mathds{R}^k)$, respectively.
From (\ref{descinex2}), we see
that $\hat L$ is differentiable and the maximal isotropy of $L$ in $\mathbf{T}^{big}_{h}M$ implies
the maximal isotropy of $\hat L$ in $\mathbf{T}^{big}_{h+k}M$. Then,
we will check the closure of $\hat L$ by the untwisted bracket
(\ref{crosettwbig}) of index $h+k$ for the three possible
combinations of components of (\ref{descinex2}). Of course, we shall
use the hypotheses on $V$, which mean that, $\forall p,q=1,...,k$,
$[V_p,V_q]=0$ and $\forall (X,u)\oplus(\alpha,v)\in\Gamma L$ one has
$(L_{V_p}X,V_pu)\oplus (L_{V_p}\alpha,V_pv)\in\Gamma L$, i.e., the
latter cross section is $G$-orthogonal to any $(X',u')\oplus(\alpha',v')\in L$. On
the other hand, we will see $\mathbf{T}^{big}_{h}M$ as a subbundle of
$\mathbf{T}^{big}_{h+k}M$ by
$$(X,u)\oplus(\alpha,v)\mapsto(X,u,0)\oplus(\alpha,v,0).$$

First we have $$ [(X_1,u_1,\alpha_1(V))\oplus(\alpha_1,v_1,0),
(X_2,u_2,\alpha_2(V))\oplus(\alpha_2,v_2,0)]_{h+k}$$ $$=
[(X_1,u_1)\oplus (\alpha_1,v_1),(X_2,u_2)\oplus (\alpha_2,v_2)]_h
\vspace*{1mm}$$
$$+(0,0,X_1(\alpha_2(V))-X_2(\alpha_1(V)))\oplus(0,0,0),$$
which has the form of a first component (\ref{descinex2}). Indeed,
after reductions, we get
$$<L_{X_1}\alpha_2-L_{X_2}\alpha_1 + \frac{1}{2}d(\alpha_1(X_2) -
\alpha_2(X_1))+\frac{1}{2}(u_2.dv_1+v_2\cdot du_1-u_1\cdot
dv_2-v_1\cdot du_2),V_p> $$ $$=X_1(\alpha_2(V_p))-X_2(\alpha_1(V_p))
+
G((L_{V_p}X_1,V_pu_1)\oplus(L_{V_p}\alpha_1,V_pv_1),(X_2,u_2)(\oplus\alpha_2,v_2))
\vspace*{1mm}$$
$$-G((L_{V_p}X_2,V_pu_2)\oplus(L_{V_p}\alpha_2,V_pv_2),(X_1,u_1)\oplus(\alpha_1,v_1))
=X_1(\alpha_2(V_p))-X_2(\alpha_1(V_p)),$$ because the vector
fields $V_p$ are infinitesimal automorphisms of $L$.

Then,
$$ [(-w_1\cdot V,0,0)\oplus(0,0,w_1),
(-w_2\cdot V,0,0)\oplus(0,0,w_2)]$$ $$ = ([w_1\cdot V,w_2\cdot V],0,0)
\oplus (0,0,w_2\cdot V(w_1)-w_1\cdot V(w_2))$$ and the commutation
of the vector fields $V_p$ shows that this result has the form of a
second component (\ref{descinex2}).

Finally, $$ [(X,u,\alpha(V))\oplus(\alpha,v,0) , (-w\cdot
V,0,0)\oplus(0,0,w)]$$ $$= (w\cdot[V,X] - (Xw)\cdot V,(w\cdot
V)(u), (w\cdot V)(\alpha(V)))\vspace*{1mm}$$ $$\oplus(w\cdot
(L_{V}\alpha),(w\cdot V)(v),Xw),$$ which is the cross
section of $\hat L$ obtained by using (\ref{Linex2}) with
$$\sum_{p=1}^kw^p(([V_p,X],V_pu)\oplus(L_{V_p}\alpha,V_pv))\in\Gamma
L$$ in the role of $(X,u)\oplus(\alpha,v)$ and $Xw$ in the role of
$w$. \end{proof}
\begin{rem}\label{obsptprelungire} {\rm
The conclusion remains true for a twisted structure $L$ if we ask
the following conditions for the vector fields $V_p$ and the twist
form (\ref{twistform}):
\begin{equation}\label{condVPhi} i(V_p)\Phi=0,\;i(V_p)\Psi=0.
\end{equation} The twist form for $\hat L$ is the same as for $L$
with the addition of $k$ forms $\Psi^p=0$.}
\end{rem}
\begin{rem}\label{prelungcombinat} {\rm Assume that $L\subseteq
\mathbf{T}^{big}_hM, L'\subseteq
\mathbf{T}^{big}_{h'}M$ are stably equivalent Dirac structures with prolongations $\hat L_1\subseteq
\mathbf{T}^{big}_{h+u}M, \hat L'_1\subseteq
\mathbf{T}^{big}_{h'+u'}M$ $(h+u=h'+u')$ isomorphic by the metric preserving automorphism $\varphi$ of $\mathbf{T}^{big}_{h+u}M$. Assume also that the commuting vector fields $V_p$ $(p=1,...,k)$ are infinitesimal automorphisms for both $L$ and $L'$. Then the prolongations 
$\hat L\subseteq
\mathbf{T}^{big}_{h+u}M, \hat L'\subseteq
\mathbf{T}^{big}_{h'+u'}M$ given by Theorem \ref{propex} are stably equivalent again. Indeed, it is easy to see that $V_p$ also are infinitesimal automorphisms of $\hat L_1,\hat L'_1$, hence, Theorem 
\ref{propex} may be applied to these structures and one gets prolongations $\widehat{(\hat L_1)},\widehat{(\hat L'_1)}\subseteq 
\mathbf{T}^{big}_{h+u+k}M$. The latter are isomorphic by the metric preserving automorphism $\hat\varphi$ of $\mathbf{T}^{big}_{h+u+k}M=
\mathbf{T}^{big}_{h+u}M\times\mathds{R}^{2k}$ that acts by $\varphi$ on the $\mathbf{T}^{big}_{h+u}M$-component and by the identity on the $\mathds{R}^{2k}$-component.}\end{rem}

We also want to add a remark of a different nature:
\begin{rem}\label{obscoh} {\rm Let $L$
be a stable Dirac structure of index $h$ and $\tilde L$ the
corresponding, translation invariant Dirac structure of
$M\times\mathds{R}^h$. A known result about Dirac structures 
tells that $d_{\tilde L}\tilde\Omega(\tilde{\mathcal{X}}_1,
\tilde{\mathcal{X}}_2,\tilde{\mathcal{X}}_3)$, where $\tilde\Omega$ is the form 
(\ref{omegabig}) on $M\times\mathds{R}^h$, $\tilde{\mathcal{X}}_a\in\Gamma\tilde L$ $(a=1,2,3)$ and $d_{\tilde L}$ is the exterior differential of the Lie algebroid $\tilde L$, is given by the twist form computed on $pr_{T(M\times\mathds{R}^h)}\tilde{\mathcal{X}}_a$. If we use left invariant arguments and the twist form (\ref{twistform}) we get
$$(d_L\Omega)( \mathcal{X}_1,
\mathcal{X}_2,\mathcal{X}_3) =\Phi(X_1,X_2,X_3) -
\sum_{Cycl(1,2,3)}u_1\cdot\Psi(X_2,X_3). $$
Thus, in the untwisted case, $\Omega|_L$ is
$d_L$-closed and there exists a {\it fundamental Lie algebroid
cohomology class} $[\Omega]\in H^2(L)$.}
\end{rem}

Now, we proceed to structures that satisfy the Dirac condition modulo a conformal change and we begin by generalities first discussed in \cite{Wd2}.
An almost Dirac structure $D\subseteq T^{big}M$ such
that $\mathcal{C}_{\tau}(D)$ is a Dirac structure, for some
$\tau\in C^\infty(M)$, i.e., $\Gamma D$ is closed by the
conformal-Courant bracket (\ref{confgen}) will be called a {\it
conformal-Dirac structure}\footnote{This name was used with a
different meaning in \cite{Wd}; on the other hand, the conformal-Dirac structures as defined here are the globally conformal Dirac structures of \cite{Wd2}.}. Obviously, $D$ and $ \mathcal{C}_{\tau}(D)$
simultaneously are maximally isotropic. Isotropy and Remark
\ref{nonC} imply that any conformal-Dirac structure is a Lie
algebroid of anchor $pr_{TM}$ with respect to the restriction of the
conformal-Dirac bracket.

Formula (\ref{conf2}) shows that an untwisted, almost Dirac
structure $D$ is conformal-Dirac iff $\forall
X\oplus\alpha,Y\oplus\beta\in\Gamma D$ one has
\begin{equation}\label{condconfD2}
[X,Y]\oplus (L_X\beta-L_Y\alpha -
d(\Omega(X\oplus\alpha,Y\oplus\beta))\end{equation}
$$+(X\tau)\beta - (Y\tau)\alpha -
(\Omega(X\oplus\alpha,Y\oplus\beta))d\tau)\in\Gamma D.$$

Let us look at the equivalent pair $(pr_{TM}(D),\theta_D)$ of the
almost Dirac structure $D$
\cite{{CW},{C}}, where
\begin{equation}\label{formapreDirac} \theta_D(X,Y)=\Omega(X\oplus
\alpha,Y\oplus\beta)=\beta(X)=-\alpha(Y)\end{equation}
$$(X,Y\in pr_{TM}(D), X\oplus
\alpha,Y\oplus\beta\in D).$$ If $D$ is a conformal-Dirac
structure then $pr_{TM}(D)=pr_{TM}(\mathcal{C}_{\tau}(D))$ is a
generalized foliation. Furthermore, (\ref{condconfD2}) means that
the $G$-scalar product of its left hand side by any $Z\oplus\gamma
\in\Gamma D$ is zero and this is easily seen to be equivalent to Wade's condition \cite{Wd2}
\begin{equation}\label{confinchis} d\theta_D =-
(d\tau|_{pr_{TM}D})\wedge\theta_D.\end{equation} Thus, a
conformal-Dirac structure produces a generalized foliation with
leaves that are {\it conformal-presymplectic manifolds} \cite{Wd2}.

More generally, if $D$ is an almost Dirac structure such that
$pr_{TM}(D)$ is a generalized foliation and the corresponding
$2$-form $\theta_D$ satisfies the condition
\begin{equation}\label{lcD} d\theta_D=\phi\wedge\theta_D,
\end{equation} where $\phi$ is a leaf-wise closed $1$-form, $D$ is a {\it locally conformal-Dirac structure} and the leaves of the corresponding generalized foliation are locally conformal presymplectic manifolds \cite{Wd2}.
\begin{example}\label{lcPoisson} {\rm\cite{Wd2} Recall the following 
formula of Gelfand and Dorfman \cite{GD}:
$$[P,P](\alpha,\beta,\gamma)=2\gamma(\sharp_P\{\alpha,\beta\}_P-
[\sharp_P\alpha,\sharp_P\beta]),$$
where $P\in\chi^2(M)$, $[P,P]$ is the Schouten-Nijenhuis bracket and
\begin{equation}\label{croset1forme} \{\alpha,\beta\}_P=
L_{\sharp_P\alpha}\beta-L_{\sharp_P\beta}\alpha
-d(P(\alpha,\beta)).\end{equation}
Accordingly, from (\ref{crosetCSW}), we get
\begin{equation}\label{formulaGD} [(\sharp_P\alpha)\oplus\alpha,
(\sharp_P\beta)\oplus\beta]\end{equation} $$=
(\sharp_P\{\alpha,\beta\}_P -
\frac{1}{2}i(\alpha\wedge\beta)[P,P])\oplus \{\alpha,\beta\}_P.$$
Now, if we use (\ref{conf2})
for $X=\sharp_P\alpha,Y=\sharp_P\beta,\Phi=0$ and ask the result
to belong to $graph\,\sharp_P$, we see that $graph\,\sharp_P$ is a 
conformal-Dirac structure iff there exists a function
$f=-2\tau\in C^\infty(M)$ such that
\begin{equation}\label{confP} [P,P]=(\sharp_Pdf)\wedge
P.\end{equation} If this happens we call $P$ a {\it conformal Poisson bivector field} or {\it structure}.
Accordingly, a {\it locally conformal Poisson structure}
is a pair $(P,\phi)$ where $P$ is a bivector
field, $\phi$ is a closed $1$-form and the following condition
holds
\begin{equation}\label{locconfP} [P,P]=(\sharp_P\phi)\wedge P.
\end{equation} For a locally conformal Poisson structure $P$,
$im\,\sharp_P$ defines a generalized foliation by locally conformal
symplectic manifolds. Hence, a locally conformal Poisson manifold is
a Jacobi manifold $(M,P,E)$, $P\in\chi^2(M),E\in\chi^1(M)$, i.e.,
\cite{DLM}
\begin{equation}\label{condJ}[P,P]=2E\wedge P,L_EP=0, \end{equation}
where $E=(1/2)\sharp_P\phi$, $d\phi=0$. The condition $L_EP=0$ is
implied by (\ref{locconfP}) because of the formula
\cite{V04}}
\begin{equation}\label{eqdinV04}
[P,P](\alpha,\beta,\phi) =2[d\phi(\sharp_P\alpha,\sharp_P\beta)-
(L_{\sharp_P\phi}P)(\alpha,\beta)].\end{equation}\end{example}

Now, we extend the notion of a Dirac-Jacobi structure defined
in \cite{{Wd},{GM}} as follows.
\begin{defin}\label{defDJ} {\rm
A {\it stable Dirac-Jacobi structure of index $h$} is a
maximally isotropic subbundle $J\subseteq\mathbf{T}^{big}_hM$
which is {\it integrable} in the sense that $J$ is closed by the
bracket (\ref{crosetwtbig}), where $\tau=\sum_{a=1}^hc_at^a$ is a linear function on $\mathds{R}^h$.}\end{defin}

A stable Dirac-Jacobi structure $J$ of index $h$ can be extended to stable Dirac-Jacobi
structures $\hat J_1,\hat J_2$ of index $h+k$ $(k\geq1)$ defined
by formula (\ref{trivialext2}) and these prolongations provide an
equivalence relation that justifies the ``stable" terminology. The structures $\hat J_1,\hat J_2$ simultaneously are integrable or not if we always use the same function $\tau$ and extend the $\Psi$ part of the twist form by $k$ zero components.

Using (\ref{deftildeL}), we can identify $J$ with a translation
invariant, almost Dirac structure $\tilde{J}$ on
$M\times\mathds{R}^h$ and get the following result.
\begin{prop}\label{integrabDJ} The untwisted, stable, almost
Dirac structure $J$ is Dirac-Jacobi iff the structure $
\mathcal{C}_\tau(\tilde{J})$ is a Dirac
structure on the manifold $M\times\mathds{R}^h$.\end{prop}
\begin{proof} The integrability of $\mathcal{C}_\tau(\tilde{J})$
obviously implies that of $J$. For the converse result, we
consider a local basis of $\mathcal{C}_\tau(\tilde{J})$ consisting
of elements of the form $\mathcal{C}_\tau((X+\sum_{a=1}^hu^a(\partial/\partial t^a))\oplus(\alpha+\sum_{a=1}^hv_adt^a))$ where $(X,u)\oplus(\alpha,v)\in J$. Then, since $\tau$ is linear, Remark
\ref{obsuzero} shows that, if $J$ is integrable, the Courant brackets of elements of this
basis belong to $\mathcal{C}_\tau(\tilde{J})$. For arbitrary
Courant brackets we get the same result if we express the
arguments of the bracket by means of the previous basis and use
property iii), Remark \ref{nonC}.\end{proof}
\begin{rem}\label{inficonform} {\rm The structure
$\mathcal{C}_\tau(\tilde{J})$ is not translation invariant.
Instead, the vector fields $\partial/\partial t^a$ are
infinitesimal conformal automorphisms of $\mathcal{C}_\tau(\tilde{J})$ in the sense of the following definition, which is interesting in its own
right. Let $D$ be an almost Dirac structure of an arbitrary
manifold $M$. A vector field $Z\in\chi^1(M)$ is an {\it
infinitesimal conformal automorphisms} of $D$ if there exists a
function $f\in C^\infty(M)$ such that
\begin{equation}\label{eqinficonf} X\oplus\alpha\in\Gamma D
\Rightarrow ([Z,X]-fX)\oplus(L_Z\alpha)\in\Gamma D.\end{equation}
This definition is motivated by the fact that if
$D=graph\,\sharp_P$ $(P\in\chi^2(M))$ then $Z$ satisfies
(\ref{eqinficonf}) iff $L_ZP=fP$.} \end{rem}
\begin{example}\label{exdeJD} {\rm
Take again the bivector field (\ref{Pi}) on $M\times\mathds{R}^h$.
If the translation invariant, almost Dirac structure
$graph\,\sharp_\Pi$ is conformal-Dirac for $\tau= \sum_{a=1}^hc_at^a$
then, by Proposition \ref{integrabDJ}, the stable almost Dirac structure $J$ such that 
$graph\,\sharp_\Pi=\tilde J$ is a stable Dirac-Jacobi
structure. The conformal-Dirac condition is equivalent with
(\ref{confP}), for $P=\Pi$ and $f=-2\sum_{a=1}^hc_at^a$, which holds
iff: i) the vector fields $V_a$ commute, ii) $[W,W]=2E\wedge W$,
iii) $[V_a,W]= E\wedge V_a$, $\forall a=1,...,h$, where
$E=\sum_{a=1}^hc_aV_a$. Conditions ii) and iii) show that $(W,E)$ is
a Jacobi structure.}\end{example}

Theorem \ref{propex} also yields a prolongation property of
Jacobi-Dirac structure.
\begin{prop}\label{propexJ}
Let $J\subseteq \mathbf{T}_h^{big}M$ be an untwisted, stable
Dirac-Jacobi structure of index $h$ on $M$ and let $V_p$
$(p=1,...,k)$ be commuting vector fields on $M$ such that
\begin{equation}\label{condptVJ} (X,u)\oplus(\alpha,v)\in\Gamma J
\end{equation} $$\Rightarrow\hspace{2mm} ([V_p,X]+d\tau(u)V_p,V_p(u))\oplus
(L_{V_p}\alpha,V_p(v)-\alpha(V_p)d\tau)\in\Gamma J.$$ Then
\begin{equation}\label{LinexJ} \hat J=\{(X-w\cdot V,u,\alpha(V))
\oplus(\alpha,v,w)\end{equation}
$$\,/\,(X,u)\oplus(\alpha,v)\in L,\,w\in\mathds{R}^k\}$$ is a
stable Dirac-Jacobi structure of index $h+k$.\end{prop}
\begin{proof} The integrability of $J$ implies that
$ \mathcal{C}_\tau(\tilde{J})$ is a Dirac structure on the
manifold $M\times\mathds{R}^h$ (remember that $\tau$ is a linear function). Condition (\ref{condptVJ}) is equivalent with
the fact that the commuting vector fields $e^{-\tau} V_p$ of
$M\times\mathds{R}^h$ are infinitesimal automorphisms of $
\mathcal{C}_\tau(\tilde{J})$ (the last term $d\tau$ of (\ref{condptVJ}) is to be seen as a vector in $\mathds{R}^h$). Hence, we can use Theorem
\ref{propex} in order to prolong $ \mathcal{C}_\tau(\tilde{J})$ to a
stable Dirac structure $\widehat{\mathcal{C}_\tau(\tilde{J})}
\subseteq
\mathbf{T}^{big}_{k}(M\times\mathds{R}^h)$. Moreover, we shall use
$e^\tau w$ instead of $w$ in the formula (\ref{Linex2}) that
defines the prolongation, which is possible since
$w\in\mathds{R}^k$ was arbitrary. The resulting prolongation can
be translated to $M\times\mathds{R}^h\times\mathds{R}^k$ and the
result is the almost Dirac structure $ \mathcal{C}_\tau(\tilde{\hat J})$ of
$M\times\mathds{R}^{h+k}$ where $\hat J$ is the structure defined
by (\ref{LinexJ}). Since by Theorem \ref{propex} this prolongation is integrable we are
done. (Notice that the function that defines the conformal change is $\tau$
and it does not depend on the coordinates on the factor $
\mathds{R}^k$.)\end{proof}

If $h=1$ and $J$ is a Jacobi structure $(W,E)$ seen as
a stable Dirac-Jacobi structure of index $1$, i.e., as a structure
defined like in Example \ref{exdeJD} with one vector field
$V_1=-E$, a technical calculation shows that hypothesis
(\ref{condptVJ}) is equivalent to
\begin{equation} [V_p,W]=E\wedge V_p,\;\;[V_p,E]=0.\end{equation}
Furthermore, the prolonged structure (\ref{LinexJ}) turns out to
be the one which is defined by the translation invariant structure
$graph\,\sharp_Q$ where \begin{equation}\label{Q}
Q=W-E\wedge\frac{\partial}{\partial t}+
\sum_{p=1}^kV_p\wedge\frac{\partial}{\partial s^p}, \end{equation}
where $s^p$ are the natural coordinates of $\mathds{R}^k$.
Hence, the prolonged structure is of the type described in Example
\ref{exdeJD} with $\tau$ having  the coefficients $c_1=-1,c_2=...=c_{k+1}=0$.
\section{Generalized almost contact structures}
Generalized complex structures recently became a subject of
interest for both geometers and physicists
\cite{{H},{G},{LMTZ}}. A generalized, almost complex structure
can be defined as a complex, almost Dirac structure $L\subseteq
T^{big}_cM =T^{big}M\otimes_{
\mathds{R}}\mathds{C}$ with the property that $L\cap\bar L=0$,
where the bar denotes complex conjugation. A necessary condition
for the existence of such a structure is the even-dimensionality
of $M$. If $L$ is integrable, i.e., closed by the (twisted)
Courant bracket, $L$ is a generalized (twisted) complex structure.

If we give a similar definition in the stable case, where
$T^{big}_cM$ is replaced by $ \mathbf{T}^{big}_cM
=(T_cM\times\mathds{C}^h)\oplus (T^*_cM\times\mathds{C}^h)$, we
get the notion of a {\it generalized (almost, twisted) stable
complex structure of index $h$}, which we will denote by a
boldface, e.g., $\mathbf{L}$. Such a structure can exist iff
$dim\,M+h$ is even.

Like the stable Dirac structures of Section 2, the stable
generalized, complex structures of index $h$ of $M$ may be
identified with translation invariant, generalized, complex
structures on the manifold $M\times\mathds{R}^h$.

We can define a prolongation of a stable, generalized, almost
complex structure $\mathbf{L}\subseteq\mathbf{T}^{big}_{c,h}M$ of index $h$ to a structure of
index $h+2k$ by taking the direct sum of the corresponding,
translation invariant structure $\tilde{\mathbf{L}}$ of
$M\times\mathds{R}^h$ with a constant complex structure of $
\mathds{R}^{2k}$. If $J_0$ is the complex structure of $
\mathds{R}^{2k}$ defined by
\begin{equation}\label{Jzer0}
J_0e_p=f_p,\;J_0f_p=-e_p,\end{equation} where $(e_p,f_p)$ is the
canonical orthonormal basis of $ \mathds{R}^{2k}$, the indicated
prolongation is defined by
\begin{equation}\label{trivialextc}\hat {\mathbf{L}}=
\{(X,u,w-\sqrt{-1}J_0w)\oplus(\alpha,v,s-\sqrt{-1}J_0s)\end{equation}
$$/\,
(X,u)\oplus(\alpha,v)\in\mathbf{L},\,w,s\in\mathds{R}^{2k}\}.$$ If
we look at the expression (\ref{crosettwbig}) for cross sections
of $\hat {\mathbf{L}}$ and use the compatibility between $J_0$ and
the scalar product in $ \mathds{R}^{2k}$, we see that $\hat
{\mathbf{L}}$ is integrable iff $
\mathbf{L}$ is integrable. (In the twisted case, we add to $\Psi$
the components $\Psi^p=0$, $p=1,...,2k$.) The definition of stably
equivalent structures via prolongations, given in Section 2 for
Dirac structures shall be adapted by the use of the prolongations
(\ref{trivialextc}).
\begin{rem}\label{obsteorema} {\rm The construction (\ref{Linex2})
of Theorem \ref{propex} can be used for a stable, generalized
complex structure $\mathbf{L}$ and for $k$ commuting, complex,
infinitesimal automorphisms $V^p$ of $\mathbf{L}$ and gives a
complex Dirac structure of index $h+k$, which is not a generalized
complex structure because
$\hat{\mathbf{L}}\cap\overline{\hat{\mathbf{L}}}\neq0$. If we
change the construction by using $2k$ commuting, complex,
infinitesimal automorphisms $V^p$ and by putting
$$\hat{\mathbf{L}}
=\{(X-(w-\sqrt{-1}J_0w)\cdot V,u,\alpha(V))
\oplus(\alpha,v,w-\sqrt{-1}J_0w)$$
$$\,/\,(X,u)\oplus(\alpha,v)\in \mathbf{L},\,w\in\mathds{R}^{2k}\},$$
we get a complex isotropic subbundle
of $\mathbf{T}^{big}_{c,h+2k}M$ such that
$\hat{\mathbf{L}}\cap\overline{\hat{\mathbf{L}}}=0$. However,
$\hat{\mathbf{L}}$ is not a generalized complex structure because
it has the complex dimension $h+k$ instead of the required
$h+2k$.}
\end{rem}
\begin{rem}\label{complexWade} {\rm We could also consider a notion
of stable, generalized, complex structure for the stable Wade
bracket, i.e., a complex almost Dirac structure
$\mathbf{J}\subseteq\mathbf{T}^{big}_{c,h}M$ which is closed by
the bracket (\ref{crosetwtbig}). In the untwisted case this is
equivalent with asking $ \mathcal{C}_\tau(\tilde{\mathbf{J}})$,
where $\tau$ is like in Definition \ref{defDJ} and
$\tilde{\mathbf{J}}$ is the corresponding translation invariant
structure, to be a generalized, complex structure on the manifold
$M\times\mathds{R}^h$ (see Proposition \ref{integrabDJ}). For $h=1$ and $\tau=t$ the structures were
discussed in
\cite{PW} under the name of {\it generalized almost contact
structures}, a name that we will prefer to use differently.}\end{rem}

Our main interest will be in the representation of a stable, generalized, almost complex
structure by classical tensor fields. As in
\cite{G}, the structure $\mathbf{L}$ is equivalent with a $G$-skew-symmetric
endomorphism $\Phi$ of $\mathbf{T}^{big}M$ such that $\Phi^2=-Id$
and the integrability condition is equivalent with the annulation of
the Nijenhuis torsion where the brackets are interpreted as stable
Courant brackets.

We shall recall the following known results
\cite{{G},{Cr},{LMTZ},{IV}}. A generalized, almost complex structure $\Phi$
on $M$ is equivalent with a triple of classical tensor fields
$(A\in\Gamma(End TM),\pi\in\chi^2(M),\sigma\in\Omega^2(M))$ obtained
from the following matrix representation of $\Phi$
\begin{equation}\label{matriceaPhi} \Phi(X\oplus\alpha)
= \left(\begin{array}{cc} A&\sharp_\pi\vspace{2mm}\\
\flat_\sigma&-^t\hspace{-1pt}A\end{array}\right)
\left(\begin{array}{l} X\vspace{2mm}\\ \alpha\end{array}\right)
=(AX+\sharp_\pi\alpha)\oplus(\flat_\sigma X-\alpha\circ A)
\end{equation} (the index $t$ denotes transposition).
The condition $\Phi^2=- Id$ is equivalent with
\begin{equation}\label{2-reddezv} A^2=-Id -
\sharp_\pi\circ\flat_\sigma,\;\pi(\alpha\circ A,\beta)=\pi(\alpha,
\beta\circ A),\;\sigma(AX,Y)=\sigma(X,AY).\end{equation} The second,
respectively the third, condition (\ref{2-reddezv}), are {\it
compatibility} of $\pi$, respectively $\sigma$, with $A$.

In terms of the classical tensor fields $(A,\pi,\sigma)$, the
integrability of $\Phi$ is expressed by the following four
conditions:

\hspace*{2mm}i) the bivector field $\pi$ defines a Poisson structure
on $M$, i.e., $[\pi,\pi]=0$;

\hspace*{1mm}ii) the {\it Schouten concomitant} of the pair
$(\pi,A)$ vanishes, i.e. (e.g., \cite{IV}),
\begin{equation}\label{SchoutenR}R_{(\pi,A)}(\alpha,X) =
\sharp_\pi(L_X(\alpha\circ A)-L_{AX}\alpha)-
(L_{\sharp_\pi\alpha}A)(X)=0;\end{equation}

iii) the Nijenhuis tensor of $A$ satisfies the condition
\begin{equation}\label{Nijptintegrab} \mathcal{N}_A(X,Y) =
[AX,AY]-A[X,AY]-A[AX,Y]+A^2[X,Y]\end{equation}
$$=\sharp_\pi[i(X\wedge Y)d\sigma];$$

iv) the {\it associated $2$-form} $\sigma_A(X,Y)=\sigma(AX,Y)$
satisfies the condition
\begin{equation}\label{difsigmaA}
d\sigma_A(X,Y,Z)=\sum_{Cycl(X,Y,Z)}d\sigma(AX,Y,Z).\end{equation}

The algebraic part of the previous results extends to stable
structures while each entry of the matrix (\ref{matriceaPhi}) will
be a $(2,2)$-matrix, corresponding to the two components of the
stable tangent bundle. We will look at the case where these
$(2,2)$-matrices contain classical tensor fields only. More
exactly, we shall assume that
\begin{equation}\label{cazulstrictclasic} A
= \left(\begin{array}{cc} F&0\vspace{2mm}\\ 0&0\end{array}\right),
\hspace{2mm} \sharp_\pi
= \left(\begin{array}{cc}\sharp_P&-^t\hspace{-1pt}Z\vspace{2mm}\\
Z&0\end{array}\right), \hspace{2mm} \flat_\sigma =
\left(\begin{array}{cc} \flat_\theta&-^t\hspace{-1pt}\xi\vspace{2mm}\\
\xi&0\end{array}\right),
\end{equation} where $F\in End(TM)$, $P\in\chi^2(M),
\theta\in\Omega^2(M)$, $Z=(Z_a):T^*M\rightarrow \mathds{R}^h$ is a
sequence of $h$ vector fields and its transposition is
$^t\hspace{-1pt}Z(u)=u\cdot Z$, and $\xi=(\xi^a):
TM\rightarrow\mathds{R}^h$ is a sequence of $h$ $1$-forms while
$^t\hspace{-1pt}\xi(u)=u\cdot \xi$.

If (\ref{cazulstrictclasic}) holds, we will say that the stable,
generalized, almost complex structure {\it has strictly classical
components} and a simple calculation shows that the conditions
(\ref{2-reddezv}) are equivalent to
\begin{equation}\label{condF}
\begin{array}{l}P(\alpha\circ F,\beta)=P(\alpha,\beta\circ F),\;
\theta(FX,Y)=\theta(X,FY),\vspace{2mm}\\F(Z_a)=0,\;\xi^a\circ
F=0,\;i(Z_a)\theta=0,\;i(\xi^a)P=0,\;\xi^a(Z_b)=\delta^a_b,\vspace{2mm}\\
F^2=-Id-\sharp_P\circ\flat_\theta+\sum_{a=1}^h\xi^a\otimes
Z_a.\end{array}\end{equation}
\begin{example}\label{excontact} {\rm
If $h=1,P=0,\theta=0$ then $(F,Z,\xi)$ is just an almost contact
structure of $M$ \cite{Bl1}. Similarly, if $h\geq1,P=0,\theta=0$
then $(F,Z_a,\xi^a)$ is a globally framed $f$-structure \cite{GY}, which we prefer to call an almost contact structure of
codimension $h$.}\end{example}

Example \ref{excontact} suggests the following definition.
\begin{defin}\label{genaprcontact} {\rm A system of tensor fields
$(P,\theta,F,Z_a,\xi^a)$ that satisfy conditions (\ref{condF}) will
be called a {\it generalized, almost contact structure of
codimension $h$}. If the corresponding, generalized, stable complex
structure is integrable the generalized, almost contact structure
will be called {\it normal}.}\end{defin}
\begin{theorem}\label{integrabaprcontact} A generalized,
almost contact structure $(P,\theta,F,Z_a,\xi^a)$ of codimension
$h$ such that $-1$ is not an eigenvalue of
$(\sharp_P\circ\flat_\theta)_x$ $(\forall x\in M)$ is normal iff
\begin{equation}\label{normalitate} \begin{array}{l}
[P,P]=0,\:R_{(P,F)}=0,\vspace{2mm}\\ L_{Z_a}P=0,\:
L_{Z_a}\theta=0,\:L_{\sharp_P\alpha}\xi^a=0,\:[Z_a,Z_b]=0,
\vspace{2mm}\\
\mathcal{N}_F(X,Y)=\sharp_P(i(X\wedge Y)d\theta) -
\sum_{a=1}^h(d\xi^a(X,Y))Z_a,\vspace{2mm}\\
d\theta_F(X,Y,Z)=\sum_{Cycl(X,Y,Z)}d\theta(FX,Y,Z).\end{array}
\end{equation}\end{theorem}
\begin{proof} In order to get the integrability conditions for the structure
(\ref{cazulstrictclasic}) we write down the integrability
conditions i)-iv) recalled above for the corresponding translation
invariant structure of $M\times\mathds{R}^h$. The classical tensor
fields (\ref{matriceaPhi}) of the translation invariant structure
are
\begin{equation}\label{Fstructura} A=F, \;
\pi=P+\sum_{a=1}^hZ_a\wedge\frac{\partial}{\partial t^a},\;
\sigma=\theta+\sum_{a=1}^h\xi^a\wedge dt^a.\end{equation}

Since $\pi$ is of the form (\ref{Pi}), we know that condition i) is
equivalent with:\\

i') $P$ is a Poisson bivector field and $Z_a$ are commuting
infinitesimal automorphisms of $P$, i.e.,
$$[Z_a,Z_b]=0,[P,P]=0,L_{Z_a}P=0.$$

Then, if we compute the Schouten concomitant
$$R_{(\pi,A)}(\alpha+\sum_{a=1}^hv_adt^a,X+\sum_{a=1}^hu^a(\partial/\partial
t^a))$$ for the four possible pairs of terms of the arguments and
use (\ref{condF}), we see that integrability condition ii)
becomes\\

ii')\hspace{1cm} $R_{(P,F)}=0$, $L_{Z_a}F=0$.\\

Similarly, if we express integrability condition iii) for the same
pairs of terms as above we get\\

iii')\hspace{1cm}
$L_{Z_a}\theta=0,\;L_{Z_b}\xi^a=0\;L_{\sharp_P\alpha}\xi^a=0\hspace{2mm}
(\forall\alpha\in\Omega^1(M))$,\\
$$ \mathcal{N}_F(X,Y)=\sharp_P(i(X\wedge Y)d\theta) -
\sum_{a=1}^h(d\xi^a(X,Y))Z_a.$$

Finally, we have
$$\sigma_A(X+\sum_{a=1}^hu^a\frac{\partial}{\partial t^a},
Y+\sum_{a=1}^hu^{'a}\frac{\partial}{\partial
t^a})=\theta(FX,Y)=\theta_F(X,Y)$$ and, if we express
(\ref{difsigmaA}) for arguments like those used to get ii'), we
see that integrability condition iv) becomes\\

iv')\hspace{1cm}$(L_{FX}\xi^a)(Y) - (L_{FY}\xi^a)(X)=0$,
$$\hspace{-20mm}d\theta_F(X,Y,Z)=\sum_{Cycl(X,Y,Z)}d\theta(FX,Y,Z).$$

Thus, i')-iv') are the normality conditions; they include the
conditions (\ref{normalitate}) and the supplementary conditions
\begin{equation}\label{condsupl} L_{Z_b}\xi^a=0,\;L_{Z_a}F=0,\;
(L_{FX}\xi^a)(Y) - (L_{FY}\xi^a)(X)=0.\end{equation} We will prove
that, if $-1$ never is an eigenvalue  of $\sharp_P\circ\flat_\theta$, (\ref{condsupl}) follow from (\ref{normalitate}).

We begin by showing that (\ref{condF}), (\ref{normalitate}) imply the existence
of some nice local coordinates and tangent bases. The condition
$\xi^a(Z_b)=\delta^a_b$ contained in (\ref{condF}) shows that the
vector fields $Z_a$ are linearly independent and so are the
$1$-forms $\xi^a$. Since the vector fields $Z_a$ commute
$span\{Z_a\}$ is tangent to a foliation $\mathcal{Z}$, with
parallelizable leaves, with local leaf-wise coordinates $z^a$ such
that $Z_a=\partial/\partial z^a$ and with transversal local
coordinates $y^u$ such that the leaves of $
\mathcal{Z}$ have the local equations $y^u=const.$ $(u=1,...,dim\,M-h)$.
Furthermore, the $1$-forms $\xi^a$ must have local expressions of
the form
\begin{equation}\label{eqluixi} \xi^a=dz^a+\xi^a_udy^u
\end{equation} and provide a complementary, normal distribution
$\nu\mathcal{Z}$ $(TM=\nu\mathcal{Z}\oplus T\mathcal{Z})$ of
equations $\xi^a=0$ and with local bases
\begin{equation}\label{bazetransv} Y_u=\frac{\partial}{\partial
y^u} - \xi^a_u\frac{\partial}{\partial z^a}.\end{equation}

Now, if we take $X=Y_u,Y=Z_b$ in the condition (\ref{normalitate}) on $\mathcal{N}_F$, we get
\begin{equation}\label{FLF}
F\circ(L_{Z_b}F)(Y_u)=-\sharp_Pi(Y_u)i(Z_b)d\theta-
\sum_{a=1}^hd\xi^a(Y_u,Z_b)Z_a\end{equation}
$$=-\sharp_Pi(Y_u)L_{Z_b}\theta +
\sum_{a=1}^h(L_{Z_b}\xi^a)(Y_u)Z_a \stackrel{(\ref{normalitate})}{=}
\sum_{a=1}^h(L_{Z_b}\xi^a)(Y_u)Z_a.$$
If we evaluate $\xi^c$ on (\ref{FLF}) and
use (\ref{condF}), we get $L_{Z_b}\xi^c(Y_u)=0$. Since it is also
simple to check that $L_{Z_b}\xi^c(Z_a)=0$, it follows that
$L_{Z_b}\xi^c=0$, which is the first condition (\ref{condsupl}).

Therefore, (\ref{FLF}) reduces to $F\circ(L_{Z_b}F)(Y_u)=0$ and if
we apply $F$ to the previous equality and use the last condition
(\ref{condF}) we get
$$(Id+\sharp_P\circ\flat_\theta)((L_{Z_b}F)(Y_u)) =
\sum_{b=1}^h\xi^a((L_{Z_b}F)(Y_u))Z_a=(L_{Z_b}\xi^a)(FY_u)=0.$$
Thus, if $\sharp_P\circ\flat_\theta$ never
has the eigenvalue $-1$ on $M$, we get $(L_{Z_b}F)(Y_u)=0$. Since in
view of (\ref{condF}), (\ref{normalitate}) we also have
$(L_{Z_a}F)(Z_b)=0$, we get $L_{Z_b}F=0$, which is the second
condition (\ref{condsupl}).

Finally, we shall derive the last condition (\ref{condsupl}). This
condition is trivial for $X=Z_a,Y=Z_b$, because of $L_{Z_a}F=0$.
For $X=Y_u,Y=Z_b$ the condition reduces to
$(L_{FY_u}\xi^a)(Z_b)=0$, which holds since
$$(L_{FY_u}\xi^a)(Z_b)=\xi^a([Z_b,FY_u])=\xi^a(L_{Z_b}F(Y_u)+F[Z_b,Y_u])=0.$$
Finally, for $X=Y_u,Y=Y_v$ one has $$(L_{FY_u}\xi^a)(Y_v) -
(L_{FY_v}\xi^a)(Y_u)= \xi^a([Y_v,FY_u]+[FY_v,Y_u]).$$ The last
expression may be calculated from the condition
(\ref{normalitate}) on $\mathcal{N}_F$, which gives
$$\xi^a(\mathcal{N}_F(FY_u,Y_v)) =-d\xi^a(FY_u,Y_v),$$ equivalently,
$$\xi^a([FY_v,F^2Y_u]+[FY_u,Y_v])=0.$$ Here, we may insert the expression of $F^2$ given by
the last condition (\ref{condF}) and the result is
$$\xi^a([Y_v,FY_u]+[FY_v,Y_u])=\xi^a([\sharp_P\circ\flat_\theta Y_u,Y_v])
=-(L_{\sharp_P\circ\flat_\theta (Y_u)}\xi^a)(Y_v)=0,$$ because of the
condition $L_{\sharp_P\alpha}\xi^a=0$ contained in
(\ref{normalitate}).\end{proof}

In the classical case $P=0,\theta=0$ Theorem \ref{integrabaprcontact} was known \cite{{Bl1},{GY}}. Theorem \ref{integrabaprcontact} and some of the facts contained
in its proof lead to the following result.
\begin{theorem}\label{strproiectata} Any normal, generalized,
almost contact structure $(P,\theta,F,Z_a,\\ \xi^a)$ of codimension
$h$ projects to a generalized complex structure of the space of
leaves of the foliation $\mathcal{Z}$.\end{theorem}
\begin{proof} In view of (\ref{condF}), the expressions of
the tensor fields $P,\theta,F$ by means of the local bases
$(Y_u,Z_a)$ and cobases $(dy^u,\xi^a)$ are of the form
\begin{equation}\label{FPtheta} P=\frac{1}{2}P^{uv}Y_u\wedge
Y_v,\;\theta=\frac{1}{2}\theta_{uv}dy^u\wedge
dy^v,\;F(Z_a)=0,\;F(Y_u)=F_u^vY_v.\end{equation} Furthermore, the
normality conditions (\ref{normalitate}) show that the
coefficients $P^{uv}$, $\theta_{uv},F_u^v$ locally depend on the
coordinates $y^u$ alone. Therefore, a projected structure defined
by $P,\theta,F$ exists. For this projected structure conditions
(\ref{2-reddezv}) and i), ii), iii) follow from the combination of (\ref{condF}) with i'), ii'), iii') since the extra terms that appear in the last condition (\ref{condF}) and in the last condition iii') are terms
in $Z_a$ and have the zero projection. Thus, the projected
structure is a generalized, complex structure.\end{proof}

In Theorem \ref{strproiectata}, if we prefer not to consider a general space of leaves, we may either refer to local transversal submanifolds of $\mathcal{Z}$ or add the hypothesis that $M/\mathcal{Z}$ is a
Hausdorff manifold. In what follows, we will give a more complete
and precise result.

In the usual way of foliation theory (e.g.,
\cite{Wk}), we define a {\it transversal, generalized, (almost) complex
structure} of a foliated manifold $(M,\mathcal{F})$ to be a maximal family $\{U_i,L_i\}$, where $\{U_i\}$ is a
covering of $M$ by $\mathcal{F}$-adapted coordinate neighborhoods
and $L_i$ is a generalized, (almost) complex structure of the local
space of leaves $U_i/(\mathcal{F}|_{U_i})$, such that, $\forall
i,j$, $L_i$ and $L_j$ restrict to the same structure on the
connected components of the open submanifold $U_i\cap
U_j/(\mathcal{F}|_{U_i\cap U_j})$. The structures $L_i$ are
equivalent with a triple of tensor fields $(A_i,\sigma_i,\pi_i)$ of
tensorial type $(1,1),(0,2),(2,0)$, respectively, on $U_i/(
\mathcal{F}|_{U_i})$, which satisfy the conditions
(\ref{2-reddezv}). If $\nu\mathcal{F}$ is a normal bundle of
$\mathcal{F}$ these tensor fields glue up to corresponding global
tensor fields $(A,\sigma,\pi)$ of the bundle $\nu\mathcal{F}$,
which have a unique extension by $0$ to projectable tensor fields
of $M$. Thus, after a choice of $\nu\mathcal{F}$, we may
equivalently define a transversal, generalized (almost) complex
structure by a projectable triple $(A,\sigma,\pi)$.

Now, we shall prove
\begin{theorem}\label{legatura} A normal, generalized, almost
contact structure $(P,\theta,F,$ $Z_a,\xi^a)$ of codimension $h$
is equivalent with the following triple of objects: {\rm1)} a
foliation $\mathcal{Z}$ endowed with a parallelization that
consists of commuting vector fields $Z_a$, {\rm2)} a
$\mathcal{Z}$-transversal, generalized complex structure
associated with $(F,P,\theta)$, {\rm3)} a normal bundle
$\nu\mathcal{Z}$, which is invariant by the linear holonomy of
$\mathcal{Z}$ and by the infinitesimal transformations belonging
to $im\,\sharp_P$ $(P\in\Gamma\wedge^2\nu\mathcal{Z})$ and has an
$F$-invariant Ehresmann curvature $(F\in
End(\nu\mathcal{Z}))$.\end{theorem}
\begin{proof} The proofs of Theorem \ref{integrabaprcontact} and Theorem
\ref{strproiectata} show that a normal structure
$(F,P,\theta,$ $Z_a,\xi^a)$ yields the three required objects. The
foliation $\mathcal{Z}$ and the parallelization $(Z_a)$ were defined
there. The $\mathcal{Z}$-transversal, generalized
complex structure is given by the projections of the tensor
fields $(F,P,\theta)$ which have the expression (\ref{FPtheta}).
The normal bundle $\nu\mathcal{Z}$ has the equations $\xi^a=0$.
The invariance properties required by 3) are ensured by the
conditions $L_{Z_b}\xi^a=0$ and $L_{\sharp_P\alpha}\xi^a=0$. The
Ehresmann curvature is defined as if $\mathcal{Z}$ would be an
Ehresmann connection of a fiber bundle, i.e., by the formula
\begin{equation}\label{Ehresmann}
R_{\nu\mathcal{Z}}(X,Y)=-pr_{T\mathcal{Z}}
[pr_{\nu\mathcal{Z}}X,pr_{\nu\mathcal{Z}}Y]\end{equation} and its
invariance by $F$ means that
\begin{equation}\label{invcurv} R_{\nu\mathcal{Z}}(FX,FY)=
R_{\nu\mathcal{Z}}(X,Y).\end{equation} If we express the last
condition iii') on $X=Y_u,Y=Y_v$ and evaluate $\xi^c$ on the
result, we get
$$[FY_u,FY_v]-[Y_u,Y_v]\in\Gamma\nu\mathcal{Z},$$ which is
equivalent to (\ref{invcurv}).

Conversely, from the objects 1) and 3) we get global $1$-forms
$\xi^a$ by asking that $\xi^a\in ann\,\nu\mathcal{Z}$ and
$\xi^a(Z_b)=\delta^a_b$. Then, we have local coordinates
where (\ref{eqluixi}) and (\ref{bazetransv}) hold and we can
define tensor fields $(F,P,\theta)$ given by (\ref{FPtheta}) that
produce the given object 2) (this is what we meant by the term ``associated" used as a condition in the formulation of 2)). The algebraic conditions
(\ref{condF}) hold. This is obvious for all but the last of them and
the last follows since the fact that
$(F,P,\theta)$ define a $\mathcal{Z}$-transversal, generalized
complex structure implies the existence of coefficients
$\lambda_u^a$ such that
$$F^2+Id+\sharp_P\circ\flat_\theta(Y_u)=\sum_{a=1}^h\lambda_u^aZ_a$$
and $\xi^a(Z_b)=\delta^a_b$ shows that $\lambda_u^a=\xi^a(Y_u)$.

Now, we shall check the integrability conditions
(\ref{normalitate}). The holonomy invariance of $\nu\mathcal{Z}$
is equivalent with $[Z_a,Y_u]\in\Gamma\nu\mathcal{Z}$ and, since
by (\ref{bazetransv}) $[Z_a,Y_u]\in\Gamma T\mathcal{Z}$, we get
$[Z_a,Y_u]=0$, which is equivalent with $L_{Z_b}\xi^a=0$. These
observations and the expressions (\ref{FPtheta}) yield
$L_{Z_a}P=0,L_{Z^a}\theta=0,L_{Z^a}F=0$.

Furthermore, the Hamiltonian invariance of $\nu\mathcal{Z}$ is
equivalent with the existence of coefficients $\kappa^a_b$ such
that $L_{\sharp_P\alpha}\xi^a=
\kappa^a_b\xi^b$. If evaluated on $Y_u$, the former condition
gives $L_{\sharp_P\alpha}\xi^a(Y_u)=0$, while
$L_{\sharp_P\alpha}\xi^a(Z_b)=0$ is a consequence of the algebraic
properties. Thus, the integrability condition
$L_{\sharp_P\alpha}\xi^a=0$ holds.

The conditions $[P,P]=0, R_{(P,F)}=0$ follow by evaluating $[P,P]$
and $R_{(P,F)}$ on the bases $(Z_a,Y_u),(\xi^a,dy^u)$ while using
the formulas (\ref{FPtheta}), (\ref{eqdinV04}) and
(\ref{SchoutenR}). If the pair of arguments that we consider are
transverse to $\mathcal{Z}$ the annulation follows because of the
integrability of the transverse, generalized, complex structure.
For other types of arguments the annulation is either
straightforward or a consequence of the holonomy and the
$P$-Hamiltonian invariance of $\nu\mathcal{Z}$ (hypotheses of
3)). A similar procedure checks the last condition
(\ref{normalitate}).

Finally, the Nijenhuis tensor condition in (\ref{normalitate}) is
checked as follows. For arguments $X=Z_a,Y=Z_b$ and $X=Z_a,Y=Y_u$
the condition can be deduced as in the calculations at the end of
the proof of Theorem \ref{integrabaprcontact}. For arguments $Y=Y_u,Y=Y_v$, because of the
integrability of the transversal, generalized, complex structure
2), there exist coefficients $\lambda_{uv}^a$ such that
$$\mathcal{N}_F(Y_u,Y_v)-\sharp_P(i(Y_v)i(Y_u)d\theta)=\sum_{a=1}^h
\lambda_{uv}^aZ_a.$$ The evaluation of $\xi^a$ on the previous
equality gives $\lambda_{uv}^a=\xi^a([FY_u,FY_v])$. Then, since
the invariance of the Ehresmann curvature (\ref{invcurv}) is
equivalent with \begin{equation}\label{condsuplF}
[FX,FY]-[X,Y]\in\Gamma\nu\mathcal{Z},\hspace{5mm}\forall X,Y
\in\Gamma\nu\mathcal{Z},\end{equation} we get
$$\lambda_{uv}^a=\xi^a([Y_u,Y_v])=-d\xi^a(Y_u,Y_v).$$
This finishes the proof.\end{proof}
\begin{rem}\label{actiuneRn} {\rm The properties of $\mathcal{Z}$ required in 1) of Theorem
\ref{legatura} are equivalent with the fact that $\mathcal{Z}$ is a
foliation of $M$ by the orbits of a locally free action of the
additive group $ \mathds{R}^h$.}\end{rem}
\begin{example}\label{fibrattoric} {\rm Let $(N,F,P,\theta)$ be a
generalized, complex manifold and $M$ a principal bundle over $N$
with the structure group $ \mathds{T}^h$, the $h$-dimensional
torus. Assume that the principal bundle $M$ has a connection $\xi$
with the connection forms $\xi^a$ and the curvature
$\Xi=\{d\xi^a\}$ satisfying the following two conditions: a)
$i(\sharp_P\alpha)\Xi=0$, b) $\Xi(FX,FY)=\Xi(X,Y)$, where $P,F$
are the lifts of the corresponding tensors of $N$ to the
horizontal distribution of the connection, extended by $0$ to non horizontal arguments. Then, there
exists a corresponding, well defined, normal, generalized, almost
contact structure of codimension $h$ on $M$. Indeed, if we denote
by $Z_a$ the fundamental, vertical, vector fields defined on $M$ by the natural basis of the Lie algebra $\mathds{R}^h$ of the structure group $\mathds{T}^h$ (e.g., see \cite{KN}) we get the object 1) required by Theorem
\ref{legatura}. The horizontal distribution of the connection may
play the role of the normal bundle $\nu\mathcal{Z}$ and the
connection forms $\xi^a$ satisfy the required algebraic
conditions. Conditions a), b) and the property $i(Z_b)\Xi^a=0$,
which holds because $\Xi$ is the curvature of a connection, ensure
the fulfillment of the remaining hypotheses of Theorem
\ref{legatura} and we are done. For instance, conditions a) and b) hold for a flat connection; this gives a concrete example of a normal, generalized, almost contact structure. Conversely, if we assume that
$(P,\theta,F,Z_a,\xi^a)$ is a normal, generalized, almost contact
structure of codimension $h$ on a compact manifold $M$ and that
suitable regularity conditions hold for the foliations defined by
the vector fields $Z_a$, we should be able to get an extended {\it
Boothby-Wang fibration theorem}
\cite{Bl1} telling that $M$ is a principal torus bundle with
connection over a manifold endowed with a usual, generalized
complex structure.}\end{example} \vspace*{2mm}

{\it Acknowledgement}. Part of the work on this paper was done during the author's visit to the Bernoulli Center of the \'Ecole Polytechnique 
F\'ed\'erale de Lausanne, Switzerland, in June-July 2006. The author wishes to express his gratitude to the Bernoulli Center and, in particular, to professor Tudor Ratiu, the director of the Center, for the invitation and support. 
 %\end{center}
\hspace*{7.5cm}{\small \begin{tabular}{l} Department of
Mathematics\\ University of Haifa, Israel\\ E-mail:
vaisman@math.haifa.ac.il \end{tabular}}
\end{document}